\documentclass[11pt,bezier]{article}
\usepackage{amsmath, amssymb, amsfonts, graphicx,tikz}
\renewcommand{\baselinestretch}{1.2}
\newtheorem{prethm}{{\bf Theorem}}
\renewcommand{\theprethm}{{\arabic{prethm}}}
\newenvironment{thm}{\begin{prethm}\sl{\hspace{-0.5
               em}{\bf.}}}{\end{prethm}}

\newtheorem{prepro}[prethm]{{\bf Proposition}}
\renewcommand{\theprepro}{{\arabic{prepro}}}

\newtheorem{prelem}[prethm]{{\bf Lemma}}
\renewcommand{\theprelem}{{\arabic{prelem}}}

\newtheorem{predeff}[prethm]{{\bf Definition}}
\renewcommand{\thepredeff}{{\arabic{predeff}}}

\newtheorem{precor}[prethm]{{\bf Corollary}}
\renewcommand{\theprecor}{{\arabic{precor}}}

\newtheorem{preconj}[prethm]{{\bf Conjecture}}
\renewcommand{\thepreconj}{{\arabic{preconj}}}

\newtheorem{preremark}[prethm]{{\bf Remark}}
\renewcommand{\thepreremark}{{\arabic{preremark}}}
\newenvironment{remark}{\begin{preremark}\rm{\hspace{-0.5
               em}{\bf.}}}{\end{preremark}}

\newtheorem{preexample}[prethm]{{\bf Example}}
\renewcommand{\thepreexample}{{\arabic{preexample}}}
\newenvironment{example}{\begin{preexample}\rm{\hspace{-0.5
               em}{\bf.}}}{\end{preexample}}

\newtheorem{preproof}{{\bf\textsf{Proof.}}}
\renewcommand{\thepreproof}{}

\newcommand{\la}{\lambda}
\newcommand{\x}{{\bf x}}
\newcommand{\y}{{\bf y}}
\newcommand{\mul}{{\rm mult}}
\newcommand{\A}{{\cal A}}
\newcommand{\q}{\textbf{\textit{q}}}
\renewcommand{\thefootnote}

\title{ Prime Labeling of Ladders}

\author{ Ebrahim Ghorbani$^{\,\rm 1, 2,}$\thanks{Corresponding author} \quad \quad  Sara Kamali$^{\,\rm 2}$ \\[.4cm]
{\sl $^{\rm 1}$Department of Mathematics, K.N. Toosi University of Technology,}\\
{\sl P.O. Box 16315-1618, Tehran, Iran}\\[0.3cm]
{\sl $^{\rm 2}$School of Mathematics, Institute for Research in Fundamental
Sciences (IPM),}\\{\sl P.O. Box
19395-5746, Tehran, Iran }
}

\begin{document}
\maketitle
\footnotetext{{\em E-mail Addresses}:  {\tt e\_ghorbani@ipm.ir} (E. Ghorbani), {\tt sakamali@ipm.ir} (S. Kamali)}

\vspace{5mm}

\begin{abstract}

A prime labeling of a graph with $n$ vertices is a labeling of  its vertices with
distinct integers from $\{1, 2,\ldots , n\}$ in such a way that the labels
of any two adjacent vertices are relatively prime.
T. Varkey conjectured that ladder graphs have a prime labeling. We prove this conjecture.

\vspace{5mm}
\noindent {\bf AMS Mathematics Subject Classification\,(2010):} 05C78
\end{abstract}

\vspace{5mm}

\section{Introduction}
Let $G$ be a simple graph with $n$ vertices.
A {\em prime labeling} of $G$ is a labeling of its vertices with
distinct integers from $\{1, 2,\ldots ,n\}$ in such a way that the labels
of any two adjacent vertices are relatively prime.
The {\em coprime graph of integers} (see \cite[Section 7.4]{ps}) has the set of all integers as vertex set where two vertices
are adjacent if and only if they are relatively prime.
So, for an $n$-vertex graph, having a prime labeling is equivalent to being a subgraph of the induced subgraph on $\{1, 2,\ldots ,n\}$  by the coprime graph of integers.
Many properties of the coprime graph of integers including identifying its subgraphs were studied by Ahlswede and Khachatrian \cite{1, 2, 3}, Erd\H{o}s \cite{7, 8, 9}, Erd\H{o}s, S\'ark\"ozy, and Szemer\'edi \cite{10, 12}, Szab\'o and T\'oth \cite{30}, Erd\H{o}s and S\'ark\"ozy \cite{11, 13},
 and S\'ark\"ozy \cite{28}. Also for a survey on the known results on this subject see \cite{ps}.

The notion of  prime labeling originated with Entringer and was introduced in \cite{t}.
Entringer around 1980 conjectured that all trees have a prime labeling. Little progress was made on this conjecture until
recently where in \cite{h} it is proved that there is an integer $n_0$ such all trees with at least $n_0$ vertices have a prime labeling.
Besides that, several classes of graphs have been shown to have a prime labeling, see \cite{g} and the aforementioned papers on coprime graphs for more details.
One of the graph classes for which the existence of prime labeling is unknown are ladders.
 Letting $P_n$ denote the path graph on $n$ vertices, the Cartesian product $P_n\times P_2$
is called the {\em $n$-ladder graph}. The purpose of this paper is to prove the following theorem.
\begin{thm}\label{main} Ladders have a prime labeling.
\end{thm}
 This result confirms a conjecture of \cite{g,v}.
Partial results on this conjecture have been reported in \cite{b,g,s1,s2,v}. 

\section{Proof of Theorem~\ref{main}}
In this section we give a proof for Theorem~\ref{main}.
Let's say that a labeled square
$$
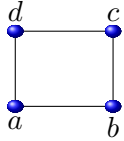
\begin{tikzpicture}[xscale=1.3]
\draw (0,0) -- (1,0);
\draw (0,1) -- (1,1);
\foreach \x in {0,1}
		\draw (\x,0) -- (\x,1);
\foreach \x in {0,1}
		\shade[ball color=blue] (\x,0) circle (.5ex);
\foreach \x in {0,1}
		\shade[ball color=blue] (\x,1) circle (.5ex);
\node [below] at (0,0) {$a$};
\node [below] at (1,0) {$b$};
\node [above] at (0,1) {$d$};
\node [above] at (1,1) {$c$};
\end{tikzpicture}$$
 fulfils the {\em cross condition} if the labels of its vertices satisfy
 $$\gcd(a,d)=\gcd(b,c)=\gcd(a,c)=\gcd(b,d)=1.$$
If we have a labeled ladder every square of which fulfils the cross condition, then we can  alternately flip the labels of vertical edges, and as a result we obtain a ladder with a prime labeling. Therefore, it suffices to prove that ladders have a labeling for which every square fulfils the cross condition.
We fix an integer $n\ge4$, and show that the $n$-ladder has such a labeling. The assertion for $n\le3$ can be verified easily.

For any positive integer $k$ we set  $$a_k:=6k-3~~~\hbox{and}~~~b_k:=6k-1.$$
We start with the following labeled $(r+s)$-ladder
$$
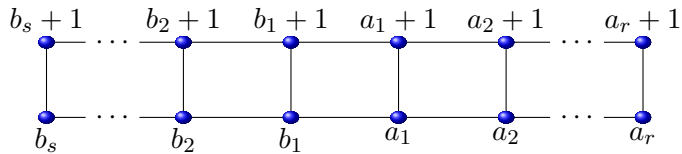
\begin{tikzpicture}[xscale=1.3]
\draw (-3.6,0) -- (-3.2,0);
\draw (-3.6,1) -- (-3.2,1);
\draw (-2.65,0) -- (1.55,0);
\draw (-2.65,1) -- (1.55,1);
\draw (2.1,0) -- (2.5,0);
\draw (2.1,1) -- (2.5,1);
\node at (-2.9,0) {$\cdots$};
\node at (-2.9,1) {$\cdots$};
\node at (1.85,0) {$\cdots$};
\node at (1.85,1) {$\cdots$};
\foreach \x in {-3.6,-2.2,-1.1,0,1.1,2.5}
		\draw (\x,0) -- (\x,1);
\foreach \x in  {-3.6,-2.2,-1.1,0,1.1,2.5}
		\shade[ball color=blue] (\x,0) circle (.5ex);
\foreach \x in  {-3.6,-2.2,-1.1,0,1.1,2.5}
		\shade[ball color=blue] (\x,1) circle (.5ex);
\node [below] at (-3.6,0) {$b_s$};
\node [below] at (-2.2,0) {$b_2$};
\node [below] at (-1.1,0) {$b_1$};
\node [below] at (0,0) {$a_1$};
\node [below] at (1.1,0) {$a_2$};
\node [below] at (2.5,0) {$a_r$};
\node [above] at (-3.6,1) {$b_s+1$};
\node [above] at (-2.2,1) {$b_2+1$};
\node [above] at (-1.1,1) {$b_1+1$};
\node [above] at (0,1) {$a_1+1$};
\node [above] at (1.1,1) {$a_2+1$};
\node [above] at (2.5,1) {$a_r+1$};
\end{tikzpicture}$$
where $r=\lfloor(n+1)/3\rfloor$ and $s=\lfloor n/3\rfloor$.
Our goal is to complete this to a labeled $n$-ladder in which every square fulfils the cross condition.

For positive integers $\ell$ satisfying
\begin{equation}\label{p1}
    (5\mid \ell+3,\,7\mid \ell+1)~~\hbox{or}~~(5\mid \ell,\,7\mid \ell+3),
\end{equation}
 we define $c_\ell:=a_\ell-2$ and we set $$C:=\{c_\ell: a_\ell-2\le2n-1\}.$$
(Note that the smallest $\ell$ satisfying \eqref{p1} is $\ell=25$, and so the smallest $n$ for which $C$ is nonempty is $73$.)
If $c_\ell$ exists, we `insert' it between $a_\ell$ and $a_{\ell+1}$ as follows.
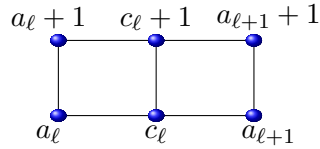
\begin{center}
\begin{tikzpicture}[xscale=1.3]
\draw (0,0) -- (2,0);
\draw (0,1) -- (2,1);
\foreach \x in {0,1,2}
		\draw (\x,0) -- (\x,1);
\foreach \x in {0,1,2}
		\shade[ball color=blue] (\x,0) circle (.5ex);
\foreach \x in {0,1,2}
		\shade[ball color=blue] (\x,1) circle (.5ex);
\node [above] at (-.1,1) {$a_\ell+1$};
\node [above] at (1,1) {$c_\ell+1$};
\node [above] at (2.15,1) {$a_{\ell+1}+1$};
\node [below] at (-.1,0) {$a_\ell$};
\node [below] at (1,0) {$c_\ell$};
\node [below] at (2.15,0) {$a_{\ell+1}$};
\end{tikzpicture}
\end{center}
In what follows by {\em inserting} an integer between two other ones, we mean the same thing as what we did above for $c_\ell$ and $a_\ell,a_{\ell+1}$.
Note that here it is possible that $a_{\ell+1}$ be larger than $2n$ and may not be presented in our labeling.
Now, let
$$D:=\{6\ell+1: 13\le6\ell+1\le2n-1\}\setminus C.$$
If $6\ell+1\in D$ and $\ell$ satisfies
\begin{equation}\label{p3}
    5\mid \ell+1~~\hbox{or}~~(5\mid \ell+3,~7\nmid \ell+1)~~\hbox{or}~~(7\mid \ell+3,~5\nmid \ell,~5\nmid \ell+2),
\end{equation}
then we insert $6\ell+1$ between $a_\ell$ and $a_{\ell+1}$:
\begin{center}
\begin{tikzpicture}[xscale=1.3]
\draw (0,0) -- (2,0);
\draw (0,1) -- (2,1);
\foreach \x in {0,1,2}
		\draw (\x,0) -- (\x,1);
\foreach \x in {0,1,2}
		\shade[ball color=blue] (\x,0) circle (.5ex);
\foreach \x in {0,1,2}
		\shade[ball color=blue] (\x,1) circle (.5ex);
\node [above] at (-.1,1) {$a_\ell+1$};
\node [above] at (1,1) {$6\ell+2$};
\node [above] at (2.15,1) {$a_{\ell+1}+1$};
\node [below] at (-.1,0) {$a_\ell$};
\node [below] at (1,0) {$6\ell+1$};
\node [below] at (2.15,0) {$a_{\ell+1}$};
\end{tikzpicture}
\end{center}
If $6\ell+1\in D$ and $\ell$ does not satisfy \eqref{p3}, then we insert $6\ell+1$ between $b_\ell$ and $b_{\ell+1}$.
 Note that if $\ell$ satisfies \eqref{p3}, it does not satisfy \eqref{p1} and so no integer from $C$ is already inserted between $a_\ell$ and $a_{\ell+1}$.
 Finally we insert $7$ between $a_1$ and $b_1$. By now, we have a labeled $(n-1)$-ladder with labels $3,\ldots,2n$.
We will insert 1  in an appropriate place later on.

Before going further with the proof, we present an example to clarify the strategy outlined above.
\begin{example} Let $n=9$ and then $r=s=3$. We first build the $(r+s)$-ladder with the following labeling:
$$
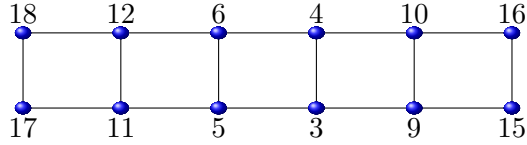
\begin{tikzpicture}[xscale=1.3]
\draw (-3,0) -- (2,0);
\draw (-3,1) -- (2,1);
\foreach \x in {-3,-2,-1,0,1,2}
		\draw (\x,0) -- (\x,1);
\foreach \x in  {-3,-2,-1,0,1,2}
		\shade[ball color=blue] (\x,0) circle (.5ex);
\foreach \x in  {-3,-2,-1,0,1,2}
		\shade[ball color=blue] (\x,1) circle (.5ex);
\node [below] at (-3,0) {$17$};
\node [below] at (-2,0) {$11$};
\node [below] at (-1,0) {$5$};
\node [below] at (0,0) {$3$};
\node [below] at (1,0) {$9$};
\node [below] at (2,0) {$15$};
\node [above] at (-3,1) {$18$};
\node [above] at (-2,1) {$12$};
\node [above] at (-1,1) {$6$};
\node [above] at (0,1) {$4$};
\node [above] at (1,1) {$10$};
\node [above] at (2,1) {$16$};
\end{tikzpicture}$$
In this $6$-ladder, the middle and far-right labeled squares do not fulfil
the cross condition. Here, the set $C$ is empty, and $D = \{13\}$. So we insert
$13$ between $a_2=9$ and $a_3=15$ yielding the following labeling of the 7-ladder:
$$
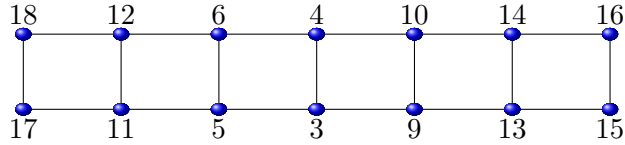
\begin{tikzpicture}[xscale=1.3]
\draw (-3,0) -- (3,0);
\draw (-3,1) -- (3,1);
\foreach \x in {-3,-2,-1,0,1,2,3}
		\draw (\x,0) -- (\x,1);
\foreach \x in  {-3,-2,-1,0,1,2,3}
		\shade[ball color=blue] (\x,0) circle (.5ex);
\foreach \x in  {-3,-2,-1,0,1,2,3}
		\shade[ball color=blue] (\x,1) circle (.5ex);
\node [below] at (-3,0) {$17$};
\node [below] at (-2,0) {$11$};
\node [below] at (-1,0) {$5$};
\node [below] at (0,0) {$3$};
\node [below] at (1,0) {$9$};
\node [below] at (2,0) {$13$};
\node [below] at (3,0) {$15$};
\node [above] at (-3,1) {$18$};
\node [above] at (-2,1) {$12$};
\node [above] at (-1,1) {$6$};
\node [above] at (0,1) {$4$};
\node [above] at (1,1) {$10$};
\node [above] at (2,1) {$14$};
\node [above] at (3,1) {$16$};
\end{tikzpicture}$$
This edge addition fixes the far-right labeled square of the former $6$-ladder so
that the cross condition now exists in the two newly formed labeled squares.
Lastly, we place a vertical edge with the label $7$
between $a_1=3$ and $b_1=5$ to create the following 8-ladder:
$$
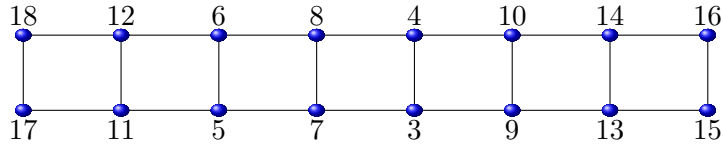
\begin{tikzpicture}[xscale=1.3]
\draw (-3,0) -- (4,0);
\draw (-3,1) -- (4,1);
\foreach \x in {-3,-2,-1,0,1,2,3,4}
		\draw (\x,0) -- (\x,1);
\foreach \x in  {-3,-2,-1,0,1,2,3,4}
		\shade[ball color=blue] (\x,0) circle (.5ex);
\foreach \x in  {-3,-2,-1,0,1,2,3,4}
		\shade[ball color=blue] (\x,1) circle (.5ex);
\node [below] at (-3,0) {$17$};
\node [below] at (-2,0) {$11$};
\node [below] at (-1,0) {$5$};
\node [below] at (0,0) {$7$};
\node [below] at (1,0) {$3$};
\node [below] at (2,0) {$9$};
\node [below] at (3,0) {$13$};
\node [below] at (4,0) {$15$};
\node [above] at (-3,1) {$18$};
\node [above] at (-2,1) {$12$};
\node [above] at (-1,1) {$6$};
\node [above] at (0,1) {$8$};
\node [above] at (1,1) {$4$};
\node [above] at (2,1) {$10$};
\node [above] at (3,1) {$14$};
\node [above] at (4,1) {$16$};
\end{tikzpicture}$$
All the squares of the 8-ladder fulfil the cross condition. Finally, we  insert 1 between any two vertices to obtain a $9$-ladder with a prime labeling.
\end{example}
\begin{remark} In what follows we frequently make use of the facts that for integers $u,v,h$,
$$\gcd(u,v)=\gcd(u,v+uh).$$
\end{remark}

Now, back to the proof, we list all the possible squares  of our labeled ladder and  verify if they fulfil the cross condition.

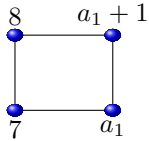
\begin{tikzpicture}[xscale=1.3]
\draw (0,0) -- (1,0);
\draw (0,1) -- (1,1);
\foreach \x in {0,1}
		\draw (\x,0) -- (\x,1);
\foreach \x in {0,1}
		\shade[ball color=blue] (\x,0) circle (.5ex);
\foreach \x in {0,1}
		\shade[ball color=blue] (\x,1) circle (.5ex);
\node [above] at (0,1) {\small$8$};
\node [above] at (1,1) {\small$a_1+1$};
\node [below] at (0,0) {\small$7$};
\node [below] at (1,0) {\small$a_1$};
\node  at (4,.5) {Clearly fulfils the cross condition.};
\end{tikzpicture}

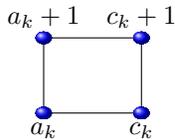
\hspace{-.4cm}\begin{tikzpicture}[xscale=1.3]
\draw (0,0) -- (1,0);
\draw (0,1) -- (1,1);
\foreach \x in {0,1}
		\draw (\x,0) -- (\x,1);
\foreach \x in {0,1}
		\shade[ball color=blue] (\x,0) circle (.5ex);
\foreach \x in {0,1}
		\shade[ball color=blue] (\x,1) circle (.5ex);
\node [above] at (0,1) {\small$a_k+1$};
\node [above] at (1,1) {\small$c_k+1$};
\node [below] at (0,0) {\small$a_k$};
\node [below] at (1,0) {\small$c_k$};
\node  at (4.2,.8) {$\gcd(a_k,c_k+1)=\gcd(a_k,a_k-1)=1$,};
\node  at (5.55,.2) {$\gcd(c_k,a_k+1)=\gcd(a_k-2,a_k+1)=\gcd(3,6k-2)=1$.};
\end{tikzpicture}

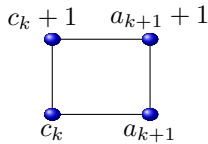
\hspace{-.4cm}\begin{tikzpicture}[xscale=1.3]
\draw (0,0) -- (1,0);
\draw (0,1) -- (1,1);
\foreach \x in {0,1}
		\draw (\x,0) -- (\x,1);
\foreach \x in {0,1}
		\shade[ball color=blue] (\x,0) circle (.5ex);
\foreach \x in {0,1}
		\shade[ball color=blue] (\x,1) circle (.5ex);
\node [above] at (-.1,1) {\small$c_k+1$};
\node [above] at (1.1,1) {\small$a_{k+1}+1$};
\node [below] at (0,0) {\small$c_k$};
\node [below] at (1,0) {\small$a_{k+1}$};
\node  at (5.7,.8) {$\gcd(c_k,a_{k+1}+1)=\gcd(6k-5,6k+4)=\gcd(9,6k+4)=1$,};
\node  at (6.85,.2) {$\gcd(a_{k+1},c_k+1)=\gcd(6k+3,6k-4)=\gcd(7,k+4)=1$ (as $k$ satisfies \eqref{p1}).};
\end{tikzpicture}

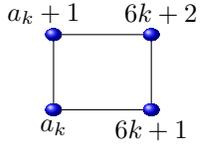
\hspace{-.4cm}\begin{tikzpicture}[xscale=1.3]
\draw (0,0) -- (1,0);
\draw (0,1) -- (1,1);
\foreach \x in {0,1}
		\draw (\x,0) -- (\x,1);
\foreach \x in {0,1}
		\shade[ball color=blue] (\x,0) circle (.5ex);
\foreach \x in {0,1}
		\shade[ball color=blue] (\x,1) circle (.5ex);
\node [above] at (-.1,1) {\small$a_k+1$};
\node [above] at (1.1,1) {\small$6k+2$};
\node [below] at (0,0) {\small$a_k$};
\node [below] at (1,0) {\small$6k+1$};
\node  at (4.8,.8) {$\gcd(6k+1,a_k+1)=\gcd(6k+1,6k+2)=1$,};
\node  at (6.8,.2) {$\gcd(a_k,6k+2)=\gcd(6k-3,6k+2)=\gcd(5,k+2)=1$ (as $k$ satisfies \eqref{p3}).};
\end{tikzpicture}

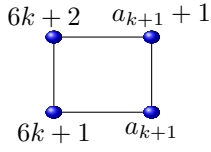
\hspace{-.4cm}\begin{tikzpicture}[xscale=1.3]
\draw (0,0) -- (1,0);
\draw (0,1) -- (1,1);
\foreach \x in {0,1}
		\draw (\x,0) -- (\x,1);
\foreach \x in {0,1}
		\shade[ball color=blue] (\x,0) circle (.5ex);
\foreach \x in {0,1}
		\shade[ball color=blue] (\x,1) circle (.5ex);
\node [above] at (-.1,1) {\small$6k+2$};
\node [above] at (1.1,1) {\small$a_{k+1}+1$};
\node [below] at (0,0) {\small$6k+1$};
\node [below] at (1,0) {\small$a_{k+1}$};
\node  at (6,.8) {$\gcd(6k+1,a_{k+1}+1)=\gcd(6k+1,6k+4)=\gcd(3,6k+4)=1$,};
\node  at (4.65,.2) {$\gcd(a_{k+1},6k+2)=\gcd(6k+3,6k+2)=1$.};
\end{tikzpicture}

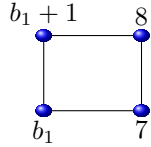
\hspace{-.2cm}\begin{tikzpicture}[xscale=1.3]
\draw (0,0) -- (1,0);
\draw (0,1) -- (1,1);
\foreach \x in {0,1}
		\draw (\x,0) -- (\x,1);
\foreach \x in {0,1}
		\shade[ball color=blue] (\x,0) circle (.5ex);
\foreach \x in {0,1}
		\shade[ball color=blue] (\x,1) circle (.5ex);
\node [above] at (0,1) {\small$b_1+1$};
\node [above] at (1,1) {\small$8$};
\node [below] at (0,0) {\small$b_1$};
\node [below] at (1,0) {\small$7$};
\node  at (4,.5) {Clearly fulfils the cross condition.};
\end{tikzpicture}

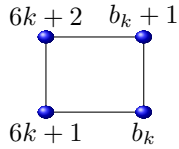
\hspace{-.2cm}\begin{tikzpicture}[xscale=1.3]
\draw (0,0) -- (1,0);
\draw (0,1) -- (1,1);
\foreach \x in {0,1}
		\draw (\x,0) -- (\x,1);
\foreach \x in {0,1}
		\shade[ball color=blue] (\x,0) circle (.5ex);
\foreach \x in {0,1}
		\shade[ball color=blue] (\x,1) circle (.5ex);
\node [above] at (0,1) {\small$6k+2$};
\node [above] at (1,1) {\small$b_k+1$};
\node [below] at (0,0) {\small$6k+1$};
\node [below] at (1,0) {\small$b_k$};
\node  at (4.5,.8) {$\gcd(6k+1,b_k+1)=\gcd(6k+1,6k)=1$,};
\node  at (4.5,.2) {$\gcd(b_k,6k+2)=\gcd(6k-1,6k+2)=1$.};
\end{tikzpicture}

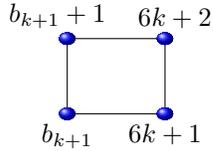
\hspace{-1.5cm}\begin{tikzpicture}[xscale=1.3]
\draw (0,0) -- (1,0);
\draw (0,1) -- (1,1);
\foreach \x in {0,1}
		\draw (\x,0) -- (\x,1);
\foreach \x in {0,1}
		\shade[ball color=blue] (\x,0) circle (.5ex);
\foreach \x in {0,1}
		\shade[ball color=blue] (\x,1) circle (.5ex);
\node [above] at (-.1,1) {\small$b_{k+1}+1$};
\node [above] at (1.1,1) {\small$6k+2$};
\node [below] at (0,0) {\small$b_{k+1}$};
\node [below] at (1,0) {{\small$6k+1$}};
\node  at (4.6,1) {$\gcd(b_{k+1},6k+2)=\gcd(6k+5,6k+2)=1$,};
\node  at (4.55,.4) {$\gcd(6k+1,b_{k+1}+1)=\gcd(6k+1,6k+6)$};
\node  at (4.55,-.2) {$\hspace{7.8cm}=\gcd(5,k+1)=1$ (as $k$ does not satisfy \eqref{p3}).};
\end{tikzpicture}

It remains to consider the following two possible squares
$$\begin{array}{ccc}
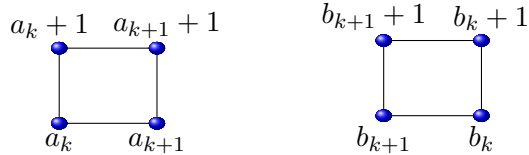
\begin{tikzpicture}[xscale=1.3]
\draw (0,0) -- (1,0);
\draw (0,1) -- (1,1);
\foreach \x in {0,1}
		\draw (\x,0) -- (\x,1);
\foreach \x in {0,1}
		\shade[ball color=blue] (\x,0) circle (.5ex);
\foreach \x in {0,1}
		\shade[ball color=blue] (\x,1) circle (.5ex);
\node [above] at (-.1,1) {$a_k+1$};
\node [above] at (1.1,1) {$a_{k+1}+1$};
\node [below] at (0,0) {$a_k$};
\node [below] at (1,0) {$a_{k+1}$};
\end{tikzpicture}
&~~~&\begin{tikzpicture}[xscale=1.3]
\draw (0,0) -- (1,0);
\draw (0,1) -- (1,1);
\foreach \x in {0,1}
		\draw (\x,0) -- (\x,1);
\foreach \x in {0,1}
		\shade[ball color=blue] (\x,0) circle (.5ex);
\foreach \x in {0,1}
		\shade[ball color=blue] (\x,1) circle (.5ex);
\node [above] at (-.1,1) {$b_{k+1}+1$};
\node [above] at (1.1,1) {$b_k+1$};
\node [below] at (0,0) {$b_{k+1}$};
\node [below] at (1,0) {$b_k$};
\end{tikzpicture}\end{array}$$
which we call  $a_k$-{\em square} and $b_k$-{\em square}, respectively.
It turns out that in some situations these two do not fulfil the cross condition, in which cases we need to replace some of the labels already assigned.
These two types of squares will be handled in the following subsections.

\subsection{$a_k$-squares}
\subsubsection{$a_k$-squares with $7\nmid k+3$ or $5\nmid k+2$}
We show that in this case,  $a_k$-squares fulfil the cross condition. Set
\begin{align*}
d_1:=\gcd(a_k,a_{k+1}+1)&=\gcd(6k-3,6k+4)=\gcd(7,k+3),\\
d_2:=\gcd(a_{k+1},a_k+1)&=\gcd(6k+3,6k-2)=\gcd(5,k+3).
\end{align*}
We have $7\nmid k+3$ or $5\nmid k+2$.
First suppose that $7\nmid k+3$. So $d_1=1$.
For a contradiction, assume that $d_2>1$, and so $5\mid k+3$.
If it happens that $6k+1=c_{k+1}\in C$, then $\ell=k+1$ satisfies \eqref{p1}, so $5\mid k+4$ or $5\mid k+1$ which is impossible.
In addition, $6k+1\ne7$. It follows that $6k+1\in D$. If $7\nmid k+1$, then $\ell=k$ satisfies \eqref{p3},  so
$6k+1$ must be already between $a_k$ and $a_{k+1}$, a contradiction. If $7\mid k+1$,  then $\ell=k$ satisfies \eqref{p1},  so
$c_k$ exists and must be already between $a_k$ and $a_{k+1}$, again a contradiction.  It follows that $d_2=1$.

Now assume that $7\mid k+3$ and $5\nmid k+2$. Again $6k+1\in D$, since
otherwise $6k+1=c_{k+1}\in C$ which means that $\ell=k+1$ satisfies \eqref{p1}, hence $7\mid k+2$ or $7\mid k+4$ which is impossible.
We also have $5\nmid k$, since otherwise $c_k$ exists, so it must had been placed between $a_k$ and $a_{k+1}$, a contradiction.
It turns out that  $\ell=k$ satisfies \eqref{p3},  so
$6k+1$ must be already between $a_k$ and $a_{k+1}$, again a contradiction.


\subsubsection{$a_k$-squares with $7\mid k+3$ and $5\mid k+2$}

In this case we have $k=35q-17$ for some positive integer $q$. Let $m=a_k=210q-105$.
It turns out that in this case $a_k$-square do not fulfil the cross condition. To overcome this obstacle, we exchange some of the labels around $a_k$.

For $k=35q-17$, it can be easily seen that $c_{k+1}$ does not exist and $6(k+1)+1=210q-95$ belongs to $D$ and $\ell=k+1$ satisfies \eqref{p3}.
  Hence $210q-95$ is already inserted between $a_{k+1}$ and $a_{k+2}$.

\noindent{\bf\textsf{Case 1.}} $11\nmid m$ and $11m>2n$.
\begin{itemize}
\item If $210q-94\le2n$, then we change the original labeling
$$
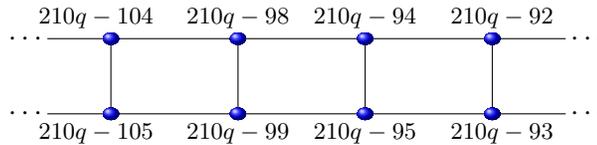
\begin{tikzpicture}[xscale=1.3]
\draw (-.65,0) -- (4.65,0);
\draw (-.65,1) -- (4.65,1);
\node at (-.85,0) {$\cdots$};
\node at (-.85,1) {$\cdots$};
\node at (4.9,0) {$\cdots$};
\node at (4.9,1) {$\cdots$};
\foreach \x in {0,1.3,2.6,3.9}
		\draw (\x,0) -- (\x,1);
\foreach \x in {0,1.3,2.6,3.9}
		\shade[ball color=blue] (\x,0) circle (.5ex);
\foreach \x in {0,1.3,2.6,3.9}
		\shade[ball color=blue] (\x,1) circle (.5ex);
\node [below] at (-.15,0) {\footnotesize$210q-105$};
\node [below] at (1.3,0) {\footnotesize$210q-99$};
\node [below] at (2.6,0) {\footnotesize$210q-95$};
\node [below] at (4,0) {\footnotesize$210q-93$};
\node [above] at (-.15,1) {\footnotesize$210q-104$};
\node [above] at (1.3,1) {\footnotesize$210q-98$};
\node [above] at (2.6,1) {\footnotesize$210q-94$};
\node [above] at (4,1) {\footnotesize$210q-92$};
\end{tikzpicture}
$$
to the following:
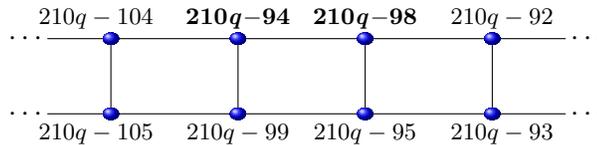
$$\begin{tikzpicture}[xscale=1.3]
\draw (-.65,0) -- (4.65,0);
\draw (-.65,1) -- (4.65,1);
\node at (-.85,0) {$\cdots$};
\node at (-.85,1) {$\cdots$};
\node at (4.9,0) {$\cdots$};
\node at (4.9,1) {$\cdots$};
\foreach \x in {0,1.3,2.6,3.9}
		\draw (\x,0) -- (\x,1);
\foreach \x in {0,1.3,2.6,3.9}
		\shade[ball color=blue] (\x,0) circle (.5ex);
\foreach \x in {0,1.3,2.6,3.9}
		\shade[ball color=blue] (\x,1) circle (.5ex);
\node [below] at (-.15,0) {\footnotesize$210q-105$};
\node [below] at (1.3,0) {\footnotesize$210q-99$};
\node [below] at (2.6,0) {\footnotesize$210q-95$};
\node [below] at (4,0) {\footnotesize$210q-93$};
\node [above] at (-.15,1) {\footnotesize$210q-104$};
\node [above] at (1.3,1) {\footnotesize{\bf210$\q-$94}};
\node [above] at (2.6,1) {\footnotesize{\bf210$\q-$98}};
\node [above] at (4,1) {\footnotesize$210q-92$};
\end{tikzpicture}$$
This changes the labeling of three consecutive squares above and as a result all three new squares  fulfil the cross condition.

\item
If $210q-94>2n$, then we change the labeling as follows:

$$\begin{array}{cc}
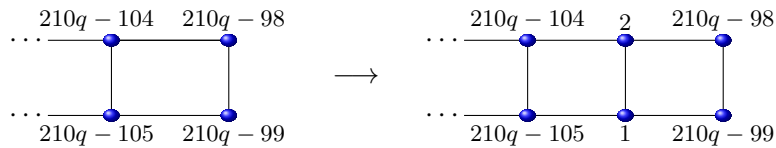
\begin{tikzpicture}[xscale=1.3]
\draw (0,0) -- (1.2,0);
\draw (0,1) -- (1.2,1);
\draw (0,1) -- (1.2,1);
\draw (-.65,1) -- (0,1);
\draw (-.65,0) -- (0,0);
\node at (-.85,0) {$\cdots$};
\node at (-.85,1) {$\cdots$};
\foreach \x in {0,1.2}
		\draw (\x,0) -- (\x,1);
\foreach \x in {0,1.2}
		\shade[ball color=blue] (\x,0) circle (.5ex);
\foreach \x in {0,1.2}
		\shade[ball color=blue] (\x,1) circle (.5ex);
\node [above] at (-.15,1) {\footnotesize$210q-104$};
\node [above] at (1.25,1) {\footnotesize$210q-98$};
\node [below] at (-.15,0) {\footnotesize$210q-105$};
\node [below] at (1.25,0) {\footnotesize$210q-99$};
\node at (2.5,.5) {$\longrightarrow$};
\end{tikzpicture}&
\begin{tikzpicture}[xscale=1.3]
\draw (-.65,0) -- (2,0);
\draw (0,1) -- (2,1);
\draw (-.65,1) -- (0,1);
\node at (-.85,0) {$\cdots$};
\node at (-.85,1) {$\cdots$};
\foreach \x in {0,1,2}
		\draw (\x,0) -- (\x,1);
\foreach \x in {0,1,2}
		\shade[ball color=blue] (\x,0) circle (.5ex);
\foreach \x in {0,1,2}
		\shade[ball color=blue] (\x,1) circle (.5ex);
\node [above] at (0,1) {\footnotesize$210q-104$};
\node [above] at (1,1) {\footnotesize$2$};
\node [above] at (2,1) {\footnotesize$210q-98$};
\node [below] at (0,0) {\footnotesize$210q-105$};
\node [below] at (1,0) {\footnotesize$1$};
\node [below] at (2,0) {\footnotesize$210q-99$};
\end{tikzpicture}
\end{array}$$

\end{itemize}

\noindent{\bf\textsf{Case 2.}} $11\nmid m$ and $11m<2n$.

Let $t$ be the maximum integer with $11^tm<2n$.

\begin{itemize}
  \item If $11\nmid m-1$, then we change the original labeling

 $$
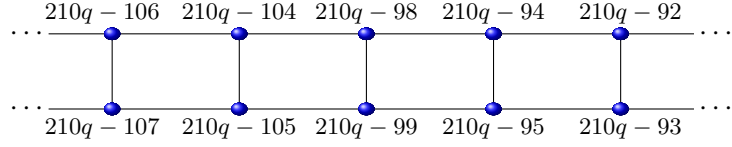
\begin{tikzpicture}[xscale=1.3]
\draw (-1.95,0) -- (4.65,0);
\draw (-1.95,1) -- (4.65,1);
\node at (-2.15,0) {$\cdots$};
\node at (-2.15,1) {$\cdots$};
\node at (4.9,0) {$\cdots$};
\node at (4.9,1) {$\cdots$};
\foreach \x in {-1.3,0,1.3,2.6,3.9}
		\draw (\x,0) -- (\x,1);
\foreach \x in {-1.3,0,1.3,2.6,3.9}
		\shade[ball color=blue] (\x,0) circle (.5ex);
\foreach \x in {-1.3,0,1.3,2.6,3.9}
		\shade[ball color=blue] (\x,1) circle (.5ex);
\node [below] at (-1.4,0) {\footnotesize$210q-107$};
\node [below] at (0,0) {\footnotesize$210q-105$};
\node [below] at (1.3,0) {\footnotesize$210q-99$};
\node [below] at (2.6,0) {\footnotesize$210q-95$};
\node [below] at (4,0) {\footnotesize$210q-93$};
\node [above] at (-1.4,1) {\footnotesize$210q-106$};
\node [above] at (0,1) {\footnotesize$210q-104$};
\node [above] at (1.3,1) {\footnotesize$210q-98$};
\node [above] at (2.6,1) {\footnotesize$210q-94$};
\node [above] at (4,1) {\footnotesize$210q-92$};
\end{tikzpicture}$$
to the following:
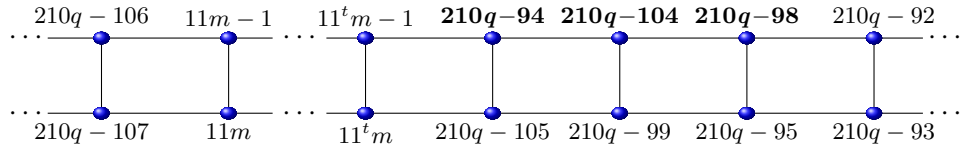
$$\begin{tikzpicture}[xscale=1.3]
\draw (-.65,0) -- (1.65,0);
\draw (2.2,0) -- (8.3,0);
\draw (-.65,1) -- (1.65,1);
\draw (2.2,1) -- (8.3,1);
\node at (-.85,0) {$\cdots$};
\node at (-.85,1) {$\cdots$};
\node at (1.95,0) {$\cdots$};
\node at (1.95,1) {$\cdots$};
\node at (8.55,0) {$\cdots$};
\node at (8.55,1) {$\cdots$};
\foreach \x in {-.1,1.2,2.6,3.9,5.2,6.5,7.8}
		\draw (\x,0) -- (\x,1);
\foreach \x in {-.1,1.2,2.6,3.9,5.2,6.5,7.8}
		\shade[ball color=blue] (\x,0) circle (.5ex);
\foreach \x in {-.1,1.2,2.6,3.9,5.2,6.5,7.8}
		\shade[ball color=blue] (\x,1) circle (.5ex);
\node [below] at (-.2,0) {\footnotesize$210q-107$};
\node [below] at (1.2,0) {\footnotesize$11m$};
\node [below] at (2.6,0) {\footnotesize$11^tm$};
\node [below] at (3.9,0) {\footnotesize$210q-105$};
\node [below] at (5.2,0) {\footnotesize$210q-99$};
\node [below] at (6.5,0) {\footnotesize$210q-95$};
\node [below] at (7.9,0) {\footnotesize$210q-93$};
\node [above] at (-.2,1) {\footnotesize$210q-106$};
\node [above] at (1.2,1) {\footnotesize$11m-1$};
\node [above] at (2.6,1) {\footnotesize$11^tm-1$};
\node [above] at (3.9,1) {\footnotesize{\bf210$\q-$94}};
\node [above] at (5.2,1) {\footnotesize{\bf210$\q-$104}};
\node [above] at (6.5,1) {\footnotesize{\bf210$\q-$98}};
\node [above] at (7.9,1) {\footnotesize$210q-92$};
\end{tikzpicture}$$

\item If $11\mid m-1$, then we change the original labeling
$$
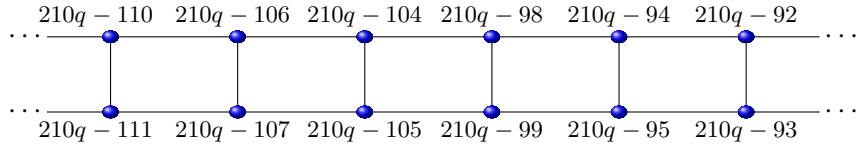
\begin{tikzpicture}[xscale=1.3]
\draw (-.65,0) -- (7.25,0);
\draw (-.65,1) -- (7.25,1);
\node at (-.85,0) {$\cdots$};
\node at (-.85,1) {$\cdots$};
\node at (7.5,0) {$\cdots$};
\node at (7.5,1) {$\cdots$};
\foreach \x in {0,1.3,2.6,3.9,5.2,6.5}
		\draw (\x,0) -- (\x,1);
\foreach \x in {0,1.3,2.6,3.9,5.2,6.5}
		\shade[ball color=blue] (\x,0) circle (.5ex);
\foreach \x in {0,1.3,2.6,3.9,5.2,6.5}
		\shade[ball color=blue] (\x,1) circle (.5ex);
\node [below] at (-.15,0) {\footnotesize$210q-111$};
\node [below] at (1.25,0) {\footnotesize$210q-107$};
\node [below] at (2.6,0) {\footnotesize$210q-105$};
\node [below] at (3.9,0) {\footnotesize$210q-99$};
\node [below] at (5.2,0) {\footnotesize$210q-95$};
\node [below] at (6.5,0) {\footnotesize$210q-93$};
\node [above] at (-.15,1) {\footnotesize$210q-110$};
\node [above] at (1.25,1) {\footnotesize$210q-106$};
\node [above] at (2.6,1) {\footnotesize$210q-104$};
\node [above] at (3.9,1) {\footnotesize$210q-98$};
\node [above] at (5.2,1) {\footnotesize$210q-94$};
\node [above] at (6.5,1) {\footnotesize$210q-92$};
\end{tikzpicture}$$
to the following:
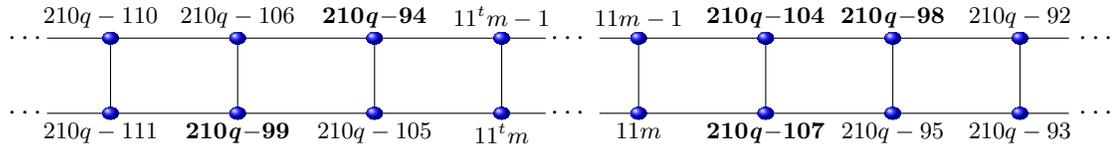
$$\begin{tikzpicture}[xscale=1.3]
\draw (-3.45,0) -- (1.65,0);
\draw (2.2,0) -- (7,0);
\draw (-3.45,1) -- (1.65,1);
\draw (2.2,1) -- (7,1);
\node at (-3.65,0) {$\cdots$};
\node at (-3.65,1) {$\cdots$};
\node at (1.9,0) {$\cdots$};
\node at (1.9,1) {$\cdots$};
\node at (7.3,0) {$\cdots$};
\node at (7.3,1) {$\cdots$};
\foreach \x in {-2.8,-1.5,-.1,1.2,2.6,3.9,5.2,6.5}
		\draw (\x,0) -- (\x,1);
\foreach \x in {-2.8,-1.5,-.1,1.2,2.6,3.9,5.2,6.5}
		\shade[ball color=blue] (\x,0) circle (.5ex);
\foreach \x in {-2.8,-1.5,-.1,1.2,2.6,3.9,5.2,6.5}
		\shade[ball color=blue] (\x,1) circle (.5ex);
\node [below] at (-2.9,0) {\footnotesize$210q-111$};
\node [below] at (-1.5,0) {\footnotesize{\bf210$\q-$99}};
\node [below] at (-.1,0) {\footnotesize$210q-105$};
\node [below] at (1.2,0) {\footnotesize$11^tm$};
\node [below] at (2.6,0) {\footnotesize$11m$};
\node [below] at (3.9,0) {\footnotesize{\bf210$\q-$107}};
\node [below] at (5.2,0) {\footnotesize$210q-95$};
\node [below] at (6.5,0) {\footnotesize$210q-93$};
\node [above] at (-2.9,1) {\footnotesize$210q-110$};
\node [above] at (-1.5,1) {\footnotesize$210q-106$};
\node [above] at (-.1,1) {\footnotesize{\bf210$\q-$94}};
\node [above] at (1.2,1) {\footnotesize$11^tm-1$};
\node [above] at (2.6,1) {\footnotesize$11m-1$};
\node [above] at (3.9,1) {\footnotesize{\bf210$\q-$104}};
\node [above] at (5.2,1) {\footnotesize{\bf210$\q-$98}};
\node [above] at (6.5,1) {\footnotesize$210q-92$};
\end{tikzpicture}$$
\end{itemize}

\noindent{\bf\textsf{Case 3.}} $11\mid m$.

As $11\mid m$, it follows that there are integers $s,q'$ such that $m=11^s(210q'-105)$ where $11\nmid 210q'-105$.
So as we did in Case 2,  $m$ and $m-1$ has  been already replaced to elsewhere. Hence we can change the labels as follows:
$$\begin{array}{cc}
\begin{tikzpicture}[xscale=1.3]
\draw (-.65,0) -- (4.65,0);
\draw (-.65,1) -- (4.65,1);
\node at (-.85,0) {$\cdots$};
\node at (-.85,1) {$\cdots$};
\node at (4.9,0) {$\cdots$};
\node at (4.9,1) {$\cdots$};
\foreach \x in {0,1.3,2.6,3.9}
		\draw (\x,0) -- (\x,1);
\foreach \x in {0,1.3,2.6,3.9}
		\shade[ball color=blue] (\x,0) circle (.5ex);
\foreach \x in {0,1.3,2.6,3.9}
		\shade[ball color=blue] (\x,1) circle (.5ex);
\node [below] at (-.15,0) {\footnotesize$210q-111$};
\node [below] at (1.25,0) {\footnotesize$210q-107$};
\node [below] at (2.65,0) {\footnotesize$210q-105$};
\node [below] at (4.05,0) {\footnotesize$210q-99$};
\node [above] at (-.15,1) {\footnotesize$210q-110$};
\node [above] at (1.25,1) {\footnotesize$210q-106$};
\node [above] at (2.65,1) {\footnotesize$210q-104$};
\node [above] at (4.05,1) {\footnotesize$210q-98$};
\node at (5.8,.5) {$\longrightarrow$};
\end{tikzpicture}
&
\begin{tikzpicture}[xscale=1.3]
\draw (-.6,0) -- (2.9,0);
\draw (-.6,1) -- (2.9,1);
\node at (-.8,0) {$\cdots$};
\node at (-.8,1) {$\cdots$};
\node at (3.15,0) {$\cdots$};
\node at (3.15,1) {$\cdots$};
\foreach \x in {0,1.15,2.3}
		\draw (\x,0) -- (\x,1);
\foreach \x in {0,1.15,2.3}
		\shade[ball color=blue] (\x,0) circle (.5ex);
\foreach \x in {0,1.15,2.3}
		\shade[ball color=blue] (\x,1) circle (.5ex);
\node [above] at (-.2,1) {\footnotesize$210q-110$};
\node [above] at (1.15,1) {\footnotesize{\bf210$\q-$104}};
\node [above] at (2.45,1) {\footnotesize$210q-98$};
\node [below] at (-.2,0) {\footnotesize$210q-111$};
\node [below] at (1.15,0) {\footnotesize{\bf210$\q-$107}};
\node [below] at (2.45,0) {\footnotesize$210q-99$};
\end{tikzpicture}
\end{array}$$

\subsection{$b_k$-squares}
\subsubsection{$b_k$-squares with $7\nmid k+1$ or $5\nmid k+1$}
We show that in this case $b_k$-squares fulfil the cross condition. Set
\begin{align*}
d_1:=\gcd(b_{k+1},b_k+1)&=\gcd(5,k),\\
d_2:=\gcd(b_k,b_{k+1}+1)&=\gcd(7,k+1).
\end{align*}
For a contradiction, suppose that  $d_1>1$,  so  $5\mid k$. We have $6k+1\in D$ since otherwise, as $6k+1\ne7$, we have $6k+1=c_{k+1}$  and so $\ell=k+1$ must satisfy \eqref{p1} which is impossible.
As $\ell=k$ does not satisfy \eqref{p3}, $6k+1$ must had been inserted between $b_k$ and $b_{k+1}$, a contradiction. Therefore, $d_1=1$.

If $7\nmid k+1$, then $d_2=1$. Assume that $7\mid k+1$ and $5\nmid k+1$.
Since $\ell=k+1$ does not satisfy \eqref{p1}, we have
necessarily $6k+1\in D$. As $5\nmid k+1$,  $\ell=k$ does not satisfy \eqref{p3}, and so $6k+1$ is already inserted between $b_k$ and $b_{k+1}$, a contradiction.

\subsubsection{$b_k$-squares with $7\mid k+1$ and $5\mid k+1$}
In this case we have $k=35q-1$ for some positive integer $q$. Let $m=b_{k+1}+1=210q$.
It turns out that in this case $b_k$-squares do not fulfil the cross condition and as before we exchange some of the labels around $b_k$.

\noindent{\bf\textsf{Case 1.}} $11\nmid m$ and $11m\ge2n$.

We change the original labeling
$$
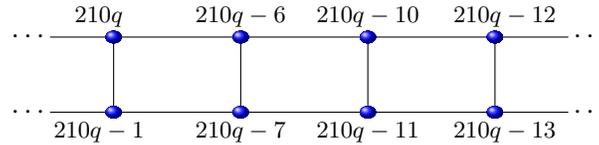
\begin{tikzpicture}[xscale=1.3]
\draw (-.65,0) -- (4.65,0);
\draw (-.65,1) -- (4.65,1);
\node at (-.85,0) {$\cdots$};
\node at (-.85,1) {$\cdots$};
\node at (4.9,0) {$\cdots$};
\node at (4.9,1) {$\cdots$};
\foreach \x in {0,1.3,2.6,3.9}
		\draw (\x,0) -- (\x,1);
\foreach \x in {0,1.3,2.6,3.9}
		\shade[ball color=blue] (\x,0) circle (.5ex);
\foreach \x in {0,1.3,2.6,3.9}
		\shade[ball color=blue] (\x,1) circle (.5ex);
\node [below] at (-.15,0) {\footnotesize$210q-1$};
\node [below] at (1.3,0) {\footnotesize$210q-7$};
\node [below] at (2.6,0) {\footnotesize$210q-11$};
\node [below] at (4,0) {\footnotesize$210q-13$};
\node [above] at (-.15,1) {\footnotesize$210q$};
\node [above] at (1.3,1) {\footnotesize$210q-6$};
\node [above] at (2.6,1) {\footnotesize$210q-10$};
\node [above] at (4,1) {\footnotesize$210q-12$};
\end{tikzpicture}
$$
to the following:
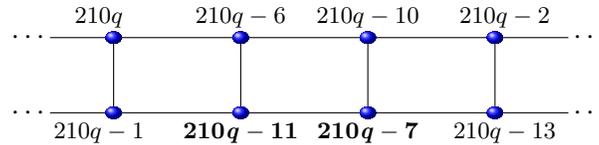
$$\begin{tikzpicture}[xscale=1.3]
\draw (-.65,0) -- (4.65,0);
\draw (-.65,1) -- (4.65,1);
\node at (-.85,0) {$\cdots$};
\node at (-.85,1) {$\cdots$};
\node at (4.9,0) {$\cdots$};
\node at (4.9,1) {$\cdots$};
\foreach \x in {0,1.3,2.6,3.9}
		\draw (\x,0) -- (\x,1);
\foreach \x in {0,1.3,2.6,3.9}
		\shade[ball color=blue] (\x,0) circle (.5ex);
\foreach \x in {0,1.3,2.6,3.9}
		\shade[ball color=blue] (\x,1) circle (.5ex);
\node [below] at (-.15,0) {\footnotesize$210q-1$};
\node [below] at (1.3,0) {\footnotesize${\bf210}\q-{\bf11}$};
\node [below] at (2.6,0) {\footnotesize${\bf210}\q-{\bf7}$};
\node [below] at (4,0) {\footnotesize$210q-13$};
\node [above] at (-.15,1) {\footnotesize$210q$};
\node [above] at (1.3,1) {\footnotesize$210q-6$};
\node [above] at (2.6,1) {\footnotesize$210q-10$};
\node [above] at (4,1) {\footnotesize$210q-2$};
\end{tikzpicture}
$$

\noindent{\bf\textsf{Case 2.}} $11\nmid m$ and $11m<2n$.

Let $t$ be the maximum integer with $11^tm<2n$.
\begin{itemize}
\item If $11\nmid m+1$, then we change the original labeling
 $$
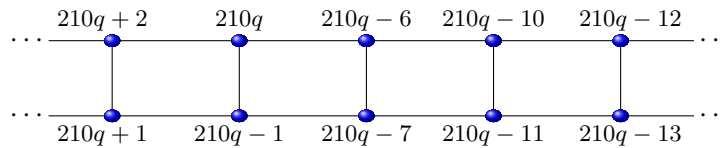
\begin{tikzpicture}[xscale=1.3]
\draw (-1.95,0) -- (4.65,0);
\draw (-1.95,1) -- (4.65,1);
\node at (-2.15,0) {$\cdots$};
\node at (-2.15,1) {$\cdots$};
\node at (4.9,0) {$\cdots$};
\node at (4.9,1) {$\cdots$};
\foreach \x in {-1.3,0,1.3,2.6,3.9}
		\draw (\x,0) -- (\x,1);
\foreach \x in {-1.3,0,1.3,2.6,3.9}
		\shade[ball color=blue] (\x,0) circle (.5ex);
\foreach \x in {-1.3,0,1.3,2.6,3.9}
		\shade[ball color=blue] (\x,1) circle (.5ex);
\node [below] at (-1.4,0) {\footnotesize$210q+1$};
\node [below] at (0,0) {\footnotesize$210q-1$};
\node [below] at (1.3,0) {\footnotesize$210q-7$};
\node [below] at (2.6,0) {\footnotesize$210q-11$};
\node [below] at (4,0) {\footnotesize$210q-13$};
\node [above] at (-1.4,1) {\footnotesize$210q+2$};
\node [above] at (0,1) {\footnotesize$210q$};
\node [above] at (1.3,1) {\footnotesize$210q-6$};
\node [above] at (2.6,1) {\footnotesize$210q-10$};
\node [above] at (4,1) {\footnotesize$210q-12$};
\end{tikzpicture}$$
to the following:
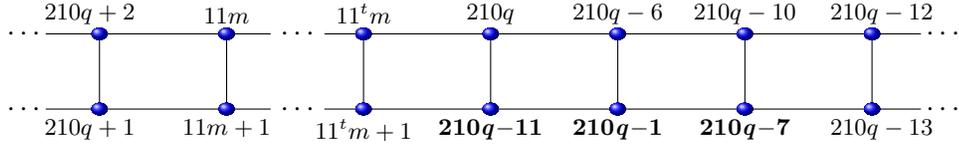
$$\begin{tikzpicture}[xscale=1.3]
\draw (-.65,0) -- (1.65,0);
\draw (2.2,0) -- (8.3,0);
\draw (-.65,1) -- (1.65,1);
\draw (2.2,1) -- (8.3,1);
\node at (-.85,0) {$\cdots$};
\node at (-.85,1) {$\cdots$};
\node at (1.95,0) {$\cdots$};
\node at (1.95,1) {$\cdots$};
\node at (8.55,0) {$\cdots$};
\node at (8.55,1) {$\cdots$};
\foreach \x in {-.1,1.2,2.6,3.9,5.2,6.5,7.8}
		\draw (\x,0) -- (\x,1);
\foreach \x in {-.1,1.2,2.6,3.9,5.2,6.5,7.8}
		\shade[ball color=blue] (\x,0) circle (.5ex);
\foreach \x in {-.1,1.2,2.6,3.9,5.2,6.5,7.8}
		\shade[ball color=blue] (\x,1) circle (.5ex);
\node [below] at (-.2,0) {\footnotesize$210q+1$};
\node [below] at (1.2,0) {\footnotesize$11m+1$};
\node [below] at (2.6,0) {\footnotesize$11^tm+1$};
\node [below] at (3.9,0) {\footnotesize{\bf210$\q-$11}};
\node [below] at (5.2,0) {\footnotesize{\bf210$\q-$1}};
\node [below] at (6.5,0) {\footnotesize{\bf210$\q-$7}};
\node [below] at (7.9,0) {\footnotesize$210q-13$};
\node [above] at (-.2,1) {\footnotesize$210q+2$};
\node [above] at (1.2,1) {\footnotesize$11m$};
\node [above] at (2.6,1) {\footnotesize$11^tm$};
\node [above] at (3.9,1) {\footnotesize$210q$};
\node [above] at (5.2,1) {\footnotesize$210q-6$};
\node [above] at (6.5,1) {\footnotesize$210q-10$};
\node [above] at (7.9,1) {\footnotesize$210q-12$};
\end{tikzpicture}$$

\item If $11\mid m+1$, then we change the original labeling
$$
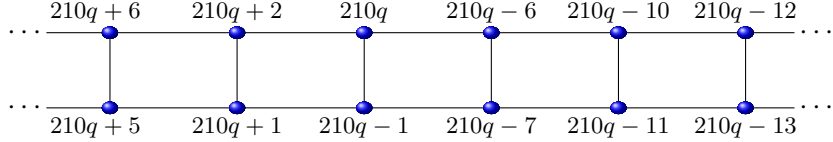
\begin{tikzpicture}[xscale=1.3]
\draw (-.65,0) -- (7,0);
\draw (-.65,1) -- (7,1);
\node at (-.85,0) {$\cdots$};
\node at (-.85,1) {$\cdots$};
\node at (7.25,0) {$\cdots$};
\node at (7.25,1) {$\cdots$};
\foreach \x in {0,1.3,2.6,3.9,5.2,6.5}
		\draw (\x,0) -- (\x,1);
\foreach \x in {0,1.3,2.6,3.9,5.2,6.5}
		\shade[ball color=blue] (\x,0) circle (.5ex);
\foreach \x in {0,1.3,2.6,3.9,5.2,6.5}
		\shade[ball color=blue] (\x,1) circle (.5ex);
\node [below] at (-.15,0) {\footnotesize$210q+5$};
\node [below] at (1.3,0) {\footnotesize$210q+1$};
\node [below] at (2.6,0) {\footnotesize$210q-1$};
\node [below] at (3.9,0) {\footnotesize$210q-7$};
\node [below] at (5.2,0) {\footnotesize$210q-11$};
\node [below] at (6.5,0) {\footnotesize$210q-13$};
\node [above] at (-.15,1) {\footnotesize$210q+6$};
\node [above] at (1.3,1) {\footnotesize$210q+2$};
\node [above] at (2.6,1) {\footnotesize$210q$};
\node [above] at (3.9,1) {\footnotesize$210q-6$};
\node [above] at (5.2,1) {\footnotesize$210q-10$};
\node [above] at (6.5,1) {\footnotesize$210q-12$};
\end{tikzpicture}$$
to the following:
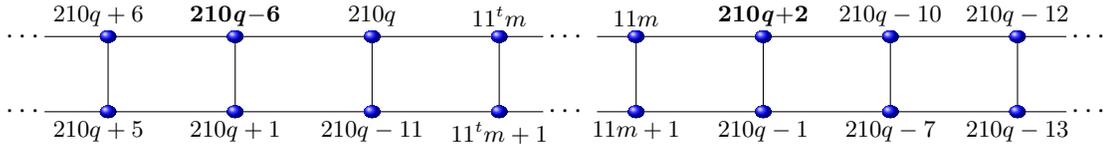
$$\begin{tikzpicture}[xscale=1.3]
\draw (-3.45,0) -- (1.65,0);
\draw (2.2,0) -- (7,0);
\draw (-3.45,1) -- (1.65,1);
\draw (2.2,1) -- (7,1);
\node at (-3.65,0) {$\cdots$};
\node at (-3.65,1) {$\cdots$};
\node at (1.9,0) {$\cdots$};
\node at (1.9,1) {$\cdots$};
\node at (7.25,0) {$\cdots$};
\node at (7.25,1) {$\cdots$};
\foreach \x in {-2.8,-1.5,-.1,1.2,2.6,3.9,5.2,6.5}
		\draw (\x,0) -- (\x,1);
\foreach \x in {-2.8,-1.5,-.1,1.2,2.6,3.9,5.2,6.5}
		\shade[ball color=blue] (\x,0) circle (.5ex);
\foreach \x in {-2.8,-1.5,-.1,1.2,2.6,3.9,5.2,6.5}
		\shade[ball color=blue] (\x,1) circle (.5ex);
\node [below] at (-2.9,0) {\footnotesize$210q+5$};
\node [below] at (-1.5,0) {\footnotesize$210q+1$};
\node [below] at (-.1,0) {\footnotesize$210q-11$};
\node [below] at (1.2,0) {\footnotesize$11^tm+1$};
\node [below] at (2.6,0) {\footnotesize$11m+1$};
\node [below] at (3.9,0) {\footnotesize$210q-1$};
\node [below] at (5.2,0) {\footnotesize$210q-7$};
\node [below] at (6.5,0) {\footnotesize$210q-13$};
\node [above] at (-2.9,1) {\footnotesize$210q+6$};
\node [above] at (-1.5,1) {\footnotesize{\bf210$\q-$6}};
\node [above] at (-.1,1) {\footnotesize$210q$};
\node [above] at (1.2,1) {\footnotesize$11^tm$};
\node [above] at (2.6,1) {\footnotesize$11m$};
\node [above] at (3.9,1) {\footnotesize{\bf210$\q+$2}};
\node [above] at (5.2,1) {\footnotesize$210q-10$};
\node [above] at (6.5,1) {\footnotesize$210q-12$};
\end{tikzpicture}$$
\end{itemize}

\noindent{\bf\textsf{Case 3.}} $11\mid m$.
\begin{itemize}
\item If $m<2n$, then $m=11^s(210q')$ for some positive integers $s,q'$ with $11\nmid 210q'$.
So as we did in Case 2,  $m$ and $m+1$ has  been already replaced to elsewhere. Hence we can change the labels as follows:
$$\begin{array}{cc}
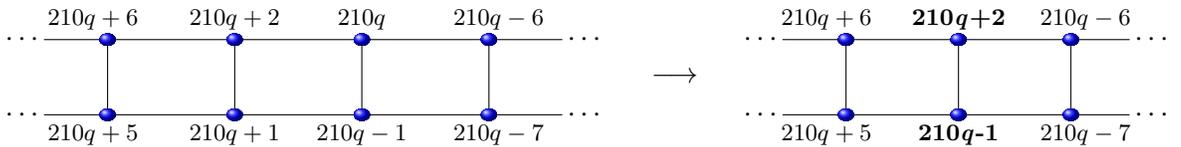
\begin{tikzpicture}[xscale=1.3]
\draw (-.65,0) -- (4.65,0);
\draw (-.65,1) -- (4.65,1);
\node at (-.85,0) {$\cdots$};
\node at (-.85,1) {$\cdots$};
\node at (4.9,0) {$\cdots$};
\node at (4.9,1) {$\cdots$};
\foreach \x in {0,1.3,2.6,3.9}
		\draw (\x,0) -- (\x,1);
\foreach \x in {0,1.3,2.6,3.9}
		\shade[ball color=blue] (\x,0) circle (.5ex);
\foreach \x in {0,1.3,2.6,3.9}
		\shade[ball color=blue] (\x,1) circle (.5ex);
\node [below] at (-.15,0) {\footnotesize$210q+5$};
\node [below] at (1.3,0) {\footnotesize$210q+1$};
\node [below] at (2.6,0) {\footnotesize$210q-1$};
\node [below] at (4,0) {\footnotesize$210q-7$};
\node [above] at (-.15,1) {\footnotesize$210q+6$};
\node [above] at (1.3,1) {\footnotesize$210q+2$};
\node [above] at (2.6,1) {\footnotesize$210q$};
\node [above] at (4,1) {\footnotesize$210q-6$};
\node at (5.8,.5) {$\longrightarrow$};
\end{tikzpicture}
&
\begin{tikzpicture}[xscale=1.3]
\draw (-.65,0) -- (2.9,0);
\draw (-.65,1) -- (2.9,1);
\node at (-.85,0) {$\cdots$};
\node at (-.85,1) {$\cdots$};
\node at (3.15,0) {$\cdots$};
\node at (3.15,1) {$\cdots$};
\foreach \x in {0,1.15,2.3}
		\draw (\x,0) -- (\x,1);
\foreach \x in {0,1.15,2.3}
		\shade[ball color=blue] (\x,0) circle (.5ex);
\foreach \x in {0,1.15,2.3}
		\shade[ball color=blue] (\x,1) circle (.5ex);
\node [above] at (-.2,1) {\footnotesize$210q+6$};
\node [above] at (1.15,1) {\footnotesize{\bf210$\q$+2}};
\node [above] at (2.45,1) {\footnotesize$210q-6$};
\node [below] at (-.2,0) {\footnotesize$210q+5$};
\node [below] at (1.15,0) {\footnotesize{\bf210$\q$-1}};
\node [below] at (2.45,0) {\footnotesize$210q-7$};
\end{tikzpicture}
\end{array}$$

\item If $m=2n$, then we insert $1$ between $b_k$ and $b_{k+1}$:
$$\begin{array}{cc}
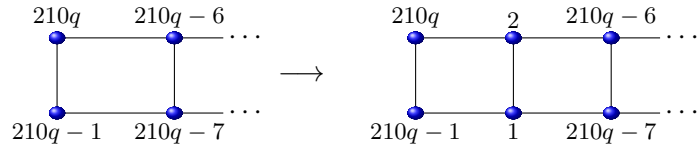
\begin{tikzpicture}[xscale=1.3]
\draw (0,0) -- (1.7,0);
\draw (0,1) -- (1.7,1);
\node at (1.95,0) {$\cdots$};
\node at (1.95,1) {$\cdots$};
\foreach \x in {0,1.2}
		\draw (\x,0) -- (\x,1);
\foreach \x in {0,1.2}
		\shade[ball color=blue] (\x,0) circle (.5ex);
\foreach \x in {0,1.2}
		\shade[ball color=blue] (\x,1) circle (.5ex);
\node [below] at (0,0) {\footnotesize$210q-1$};
\node [below] at (1.25,0) {\footnotesize$210q-7$};
\node [above] at (0,1) {\footnotesize$210q$};
\node [above] at (1.25,1) {\footnotesize$210q-6$};
\node at (2.5,.5) {$\longrightarrow$};
\end{tikzpicture}&
\begin{tikzpicture}[xscale=1.3]
\draw (0,0) -- (2.5,0);
\draw (0,1) -- (2.5,1);
\node at (2.75,0) {$\cdots$};
\node at (2.75,1) {$\cdots$};
\foreach \x in {0,1,2}
		\draw (\x,0) -- (\x,1);
\foreach \x in {0,1,2}
		\shade[ball color=blue] (\x,0) circle (.5ex);
\foreach \x in {0,1,2}
		\shade[ball color=blue] (\x,1) circle (.5ex);
\node [below] at (0,0) {\footnotesize$210q-1$};
\node [below] at (1,0) {\footnotesize$1$};
\node [below] at (2,0) {\footnotesize$210q-7$};
\node [above] at (0,1) {\footnotesize$210q$};
\node [above] at (1,1) {\footnotesize$2$};
\node [above] at (2,1) {\footnotesize$210q-6$};
\end{tikzpicture}
\end{array}$$
\end{itemize}

Finally if the labels 1 and 2 have not been used, we may insert 1 between any two vertices.

\section*{Acknowledgments}

 The research of the author was in part supported by a grant from IPM.

{}

\end{document}